\documentclass[12pt]{amsart}

\usepackage{amsfonts,amsthm,amsmath,amssymb,graphicx, overpic}

\usepackage{enumitem}

\usepackage{graphicx}

\makeatletter
\@namedef{subjclassname@2020}{%
  \textup{2020} Mathematics Subject Classification}
\makeatother

\usepackage[T1]{fontenc}

\usepackage[usenames]{color}

\newtheorem{theorem}{Theorem}[section]

\newtheorem{lemma}[theorem]{Lemma}

\theoremstyle{definition}
\newtheorem{definition}[theorem]{Definition}
\newtheorem{remark}[theorem]{Remark}

\numberwithin{equation}{section}

\newcommand{\R}{\mathbb R}

\newcommand{\N}{\mathbb N}

\renewcommand{\epsilon}{\varepsilon}

\begin{document}


\baselineskip=17pt


\title[Big Birkhoff sums in $d$-decaying Gauss like systems]{Big Birkhoff sums in $d$-decaying  Gauss like iterated function systems}
\date{}

\author[L. Liao]{Lingmin Liao}
\address{Univ Paris Est Creteil, CNRS, LAMA, F-94010 Creteil,
  France\\ Univ Gustave Eiffel, LAMA, F-77447 Marne-la-Vall\'ee,
  France}
\email{lingmin.liao@u-pec.fr}
\author[M. Rams]{Micha\l\ Rams}
\address{Institute of Mathematics\\ Polish
Academy of Sciences\\ 
\'Sniadeckich 8\\ 00-656 Warszawa, Poland }
\email{rams@impan.pl}

\date{}

\begin{abstract}
The increasing rate of the Birkhoff sums in the infinite iterated function systems with polynomial decay of the derivative (for example the Gauss map) is studied. For different unbounded potential functions, the Hausdorff dimensions of the sets of points whose Birkhoff sums share the same increasing rate are obtained. 


\end{abstract}



\def\thefootnote{}
\def\thefootnote{\arabic{footnote}}

\maketitle

\section{Introduction}


 Denote by $\N=\{1,2, \dots \}$ the set of positive integers. Consider the so-called Gauss infinite iterated function system $\left\{T_n\right\}_{n\in \mathbb{N}}$ on the unit interval $[0,1]$ defined by
\[
T_n(x):= {1 \over x+n} \ \ \text{for} \ x\in [0,1].
\]
It is well-known that the limit set of the Gauss iterated function system is the set of all irrational numbers in the unit interval $[0,1]$. In fact, for any $x\in [0,1]\setminus \mathbb{Q}$, there exists a unique infinite sequence $(a_1, a_2, \dots) \in \mathbb{N}^\mathbb{N}$ satisfying
\[
x = \lim_{n\to \infty} T_{a_1}\circ \cdots \circ T_{a_n}(1)=\frac{\displaystyle 1}{\displaystyle a_1+ \frac{\displaystyle
1}{\displaystyle a_2+\frac{\displaystyle 1}{\displaystyle
a_3+\ddots}}}.
\]
The latter is the regular continued fraction expansion of $x$. The digits $a_j=a_j(x)$ are called the partial quotients of $x$ in its continued fraction expansion. 
%

\smallskip
 For any $n\geq 1$, denote $S_n(x)=\sum_{j=1}^n a_j(x)$. In the literature, we are interested in the sum $S_n(x)$ of partial quotients which is a special Birkhoff sum with respect to the Gauss  iterated function system $\{T_n\}_{n\geq 1}$. In fact, the functions $T_n$ are inverse branches of the Gauss
transformation $T:[0,1]\to [0,1]$ defined by  $$ T(0):=0, \ \text{and } \
T(x):=\frac{1}{x} \ {\rm{(mod \ 1)}}, \ {\rm{for}}\ x\in (0,1].
$$
Then $a_1(x)= \lfloor x^{-1}\rfloor$ ($\lfloor \cdot \rfloor $ stands for the integer part)
  and $a_j(x)=a_1(T^{j-1}(x))$ for $j\geq 2$. Thus 
\[
S_n(x)= \sum_{j=1}^n a_j(x)= \sum_{j=1}^n a_1(T^{j-1}x)
\]  
is a Birkhoff sum of the potential function $x \mapsto a_1(x)$ with respect to the Gauss transformation.

In 1935, Khintchine \cite{Kh35} showed that $S_n(x)/(n\log n)$ converges in measure (Lebesgue measure) to the constant $1/\log 2$. In 1988, Philipp \cite{Ph88} proved that there is no normalizing sequence $\Phi(n)$ such that $S_n(x)/\Phi(n)$ converges to a positive constant Lebesgue almost surely. Motivated by such a phenomenon, people then turn to study the sums $S_n(x)$ from the point of view of multifractal analysis. Precisely, one is concerned with the Hausdorff dimension of the sets
\[
E(\Phi) = \left\{x\in (0,1): \lim_{n\to\infty} \frac {S_n(x)} {\Phi(n)} =1\right\},
\]
where $\Phi :\N\rightarrow \R^+$ is an increasing function. 
When $\Phi(n)/n$ has a finite limit as $n\to \infty$, $E(\Phi)$ is the classical level set of Birkhoff averages, and its Hausdorff dimension has been determined by Iommi and Jordan \cite{IJ}.  The particular attention is thus paid to the cases when the sums $S_n(x)$ are bigger, that is, when 
\[
\lim_{n\to \infty} {\Phi(n) \over n}=\infty.
\]
In this direction,  if $\Phi(n)=n^{a}$ with $a \in (1,\infty)$ or $\Phi(n)=\exp(n^\alpha)$ with $\alpha\in (0,1/2)$, Wu and Xu \cite{WX11} proved that $\dim_HE(\Phi)=1$. Here and in what follows, $\dim_H$ stands for the Hausdorff dimension. Later, Xu \cite{Xu} proved that if $\Phi(n)=\exp(n^\alpha)$ with $\alpha\in [1, \infty)$ then $\dim_HE(\Phi)=1/2$; and if $\Phi(n)=\exp(\beta^n)$ with $\beta>1$ then $\dim_HE(\Phi)=1/(\beta+1)$.  The gap for the case  $\Phi(n)=\exp(n^\alpha)$ with $\alpha\in [1/2, 1)$ was finally filled by the authors in \cite{LR16a} where we proved that $\dim_HE(\Phi)=1/2$ for all  $\alpha\in [1/2, 1)$. 
Hence there is a jump of Hausdorff dimension from $1$ to $1/2$ for the class $\Phi(n) = \exp(n^\alpha)$ at $\alpha = 1/2$. 

%
%
%
%

\smallskip
The present paper aims at generalizing the above results on the Birkhoff sums of the potential $x\mapsto a_1(x)$ in Gauss infinite iterated function system associated to continued fractions to Birkhoff sums of a general potential function in some general infinite function systems. We are especially interested in big Birkhoff sums. Before stating our main results, let us give some notations and definitions. 

\begin{definition}\label{Def:d-decaying}
 Let $d>1$ be a real number. A family $\{f_n\}_{n\in\N}$ of $C^1$ maps from the interval $[0, 1]$ to itself is called a {\it $d$-decaying  Gauss like} iterated function system  if the following properties are satisfied:
\begin{itemize}

\item[(1)] for any $i, j \in \N$, \  $f_i((0, 1)) \cap f_j((0, 1)) = \emptyset$;
\item[(2)] $\bigcup_{i=1}^\infty f_i([0,1]) = [0,1)$; 
\item[(3)]  if $f_i(x) < f_j(x)$ for all $x\in (0,1)$ then $i<j$;
\item[(4)] there exists $m\in \N$ and $0 < A < 1$ such that for all $(a_1, . . . , a_m) \in \N^m$ and for all $x \in [0, 1]$
\[
0 < |(f_{a_1} \circ \cdots \circ f_{a_m})'(x)| \leq A < 1;
\]
\item[(5)] for any $\delta > 0$, we can find two constants $K_1=K_1(\delta), K_2=K_2(\delta)>0$ such that for $i \in \N$ there exist constants $\xi_i, \lambda_i$ such that 
\[\forall x \in [0, 1], \qquad \frac {K_1} {i^{d+\delta}} \leq  \xi_i \leq |f_i'(x)| \leq \lambda_i  \leq \frac {K_2} {i^{d-\delta}}.\]
\end{itemize}
\end{definition}

We have a natural projection $\Pi: \N^\N \to [0,1]$ defined by
\[
\Pi(\underline{a}) = \lim_{n\to \infty} f_{a_1}\circ \cdots \circ f_{a_n}(1).
\]
The image $\Pi(\mathbb{N}^\mathbb{N})$ is called the limit set of the iterated function system. The set $\Pi(\N^\N)$ is always equal to $[0,1]$ minus an at most countable set.  Moreover, the map $\Pi$ is bijective, except possibly at a countable set of points where it can be $2$-to-$1$. See \cite[Chapter 1]{MU} for more details.
The inverse of $\Pi$ gives for points $x\in [0,1]$ their symbolic expansions in $\mathbb{N}^{\mathbb{N}}$. The symbolic expansion is unique for most points, but there can exist countably many points that have zero or two symbolic expansions. When the symbolic expansion is unique, we write $x=(a_1(x), a_2(x), \ldots)$ the expansion of $x\in [0,1]$. 



%
%
For each $n\in \mathbb{N}$, and each word $a_1\cdots a_n\in \N^n$, the set
\[
I_n(a_1, \cdots, a_n)=f_{a_1}\circ\cdots\circ f_{a_n}([0,1])
\]
is called a basic interval of order $n$ or an $n$-basic interval. Except for a countable set, the $n$-basic interval $I_n(a_1, \cdots, a_n)$ is identical with the set of points $x\in [0,1]$ whose symbolic expansions begin with $a_1, \cdots, a_n$. Write $I_n(x)$ the $n$-basic interval containing $x\in [0,1]$.

Denote by $|I|$ the diameter of an interval $I$. 
We say the $d$-decaying Gauss like iterated function system $\{f_n\}_{n\in \mathbb{N}}$ satisfies the {\it bounded distortion property} if there exist positive constants $K_3$ and $K_4$ such that for any two finite words $a_1a_2 \cdots a_n\in\mathbb{N}^n$ and $b_1b_2\cdots b_m\in\mathbb{N}^m$, we have
\begin{align}\label{distortion}
K_3  \leq { |I_{n+m}(a_1,\cdots, a_n,b_1,\cdots, b_m)| \over |I_n(a_1,\cdots, a_n)| \cdot |I_m(b_1,\cdots, b_m)|} \leq K_4. 
\end{align}

{ We remark that the bounded distortion property is usually stated as that there exists a constant $K_5>0$ such that for any word $a_1\cdots a_n\in\mathbb{N}^n$ and any $x,y\in I_n(a_1,\ldots,a_n)$,
\begin{align}\label{distortion-2}
{1 \over K_5}  \leq  \left| { (f_{a_1}\circ\cdots\circ f_{a_n})'(x) \over (f_{a_1}\circ\cdots\circ f_{a_n})'(y)} \right|\leq K_5. 
\end{align}
By the Mean Value Theorem, one can easily check that if (\ref{distortion-2}) holds then (\ref{distortion}) holds with $K_3=1/K_5^2$, and $K_4=K_5^2$.

The bounded distortion is satisfied for any finite or countable conformal iterated function systems, in particular it is satisfied for the classical Gauss infinite iterated function system. However, it is usually not satisfied for nonuniformly hyperbolic iterated function systems, like the iterated function systems associated to the Manneville-Pomeau map and the backward continued fraction map.

}

Consider a {\it potential function} $\varphi :[0,1] \rightarrow\R_+$, such that $\varphi$ is a constant on the interior of $I_1(j)$ for all $j\in \N$. {For all $j\in\N$, for simplicity, in what follows, we denote by $\varphi(j)$ the constant value of $\varphi$ on the interior of $1$-basic interval $I_1(j)$.  Without confusion, we also consider $\varphi$ as a function from $\mathbb{N}$ to $\R_+$.}
For $n\in \mathbb{N}$, the $n$-th Birkhoff sum of $\varphi$ at $x\in (0,1)$ is defined by
$$S_n\varphi(x)=\sum_{k=1}^{n}{ \varphi(a_k)}, \quad \text{if} \ x\in I_n(a_1, \cdots, a_n). 
$$ 
We remark that except for a countable set, the above Birkhoff sums are well defined.

For a positive {\it growth rate function} $\Phi: \mathbb{N} \rightarrow \mathbb{R}_+$, we are interested in the following set
\begin{align}
E_\varphi(\Phi):=\left\{ x\in (0,1): \lim_{n\to\infty} {S_n\varphi(x)\over \Phi(n)}=1\right\}.
\end{align}


We will calculate $\dim_HE_\varphi(\Phi)$. 
As in the Gauss iterated function system, when $\Phi(n)/n$ has a finite limit as $n\to \infty$, the set $E_\varphi(\Phi)$ is the classical level set of Birkhoff averages and its Hausdorff dimension has been well studied in \cite{PW99, KS07, FLWW, IJ, FJLR, KLRW}, and many other papers. In this paper we will consider the case when $\Phi(n)/n\to \infty$ as $n\to \infty$, thus necessarily the potential function $\varphi$ is unbounded in $[0,1]$. 


%

We obtain the following multifractal analysis results on the Hausdorff dimension of $E_\varphi(\Phi)$, according to different choices of $\varphi$ and $\Phi$. We will see that as in the case of Gauss iterated function system, the jump of Hausdorff dimension also happens but at different places. 

\begin{theorem}\label{thm-1}
Suppose $\varphi(j) = {j^a}$ for all $j\geq 1$, with $a>0$.\\
{\rm (I)} When $\Phi(n)=e^{n^\alpha}$ with $\alpha>0$, we have \\
{\rm (I-1)} $\dim_HE_\varphi(\Phi)=1$ if $\alpha < {1 \over 2}$ and the distortion property (\ref{distortion}) holds;\\
{\rm (I-2)} $\dim_HE_\varphi(\Phi)=1/d$ if $\alpha> {1 \over 2}$.\\
{\rm (II)} When  $\Phi(n)=e^{\beta^n}$ with $\beta>1$, we have $\dim_HE_\varphi(\Phi)={1\over d\beta -\beta +1}$.
\end{theorem}

\begin{theorem}\label{thm-2}
Suppose $\varphi(j) = e^{(\log j)^b}$ for all $j\geq 1$, with $b> 1$.\\
{\rm (I)} When $\Phi(n)=e^{n^\alpha}$ with $\alpha>0$, we have \\
{\rm (I-1)} $\dim_HE_\varphi(\Phi)=1$ if $\alpha < {b \over b+1}$ and the distortion property (\ref{distortion}) holds;\\
{\rm (I-2)} $\dim_HE_\varphi(\Phi)=1/d$ if $\alpha >{b \over b+1} $.\\
{\rm (II)} When $\Phi(n)={e^{\beta^n}}$ with $\beta>1$, we have $\dim_HE_\varphi(\Phi)={1\over d\beta^{1 \over b} -\beta^{1 \over b} +1}$.
\end{theorem}

\begin{theorem}\label{thm-3}
Suppose $\varphi = e^{j^c}$ for all $j\geq 1$, with $0<c<1$.\\
{\rm (I)}  When $\Phi(n)=e^{n^\alpha}$ with $\alpha>0$, we have \\
{\rm (I-1)}  $\dim_HE_\varphi(\Phi)=1$ if $\alpha < 1$ and the distortion property (\ref{distortion}) holds;\\
{\rm (I-2)}  $\dim_HE_\varphi(\Phi)={1-c \over d}$ if $\alpha> 1$. \\
{\rm (II)}  When $\Phi(n)={e^{\beta^n}}$ with $\beta>1$, we have $\dim_HE_\varphi(\Phi)={1-c \over d}$.\\
{\rm (III)}  When $\Phi(n)=e^{e^{\gamma^n}}$ with $\gamma>1$, we have $\dim_HE_\varphi(\Phi)={1-c \over d\gamma -(1-c)(\gamma -1)}$.
\end{theorem}

\begin{theorem}\label{thm-4}
Suppose $\varphi(j) = e^{j^c}$ for all $j\geq 1$, with $c\geq 1$. When $\Phi(n)=e^{n^\alpha}$, with $\alpha>0$, we have \\
{\rm (I-1)}  $\dim_HE_\varphi(\Phi)=1$ if $\alpha < 1$ and the distortion property (\ref{distortion}) holds;\\
{\rm (I-2)}  $\dim_HE_\varphi(\Phi)=0$ if $\alpha \geq 1$.
\end{theorem}

\medskip
The Hausdorff dimensions in Theorems \ref{thm-1}-\ref{thm-4} are depicted in Figures 1-4.

\begin{figure}
\vskip 25pt
\begin{overpic}[scale=1,
]{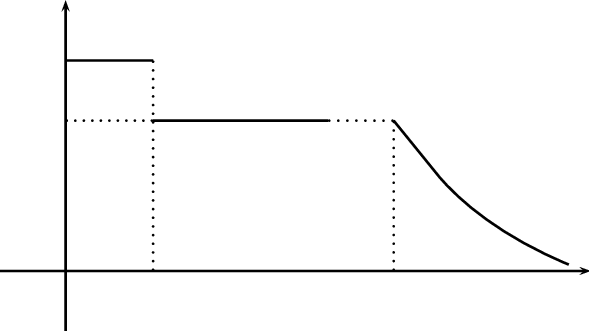}
   \put(64,5){$\beta=1+$}
 
  \put(8.8,5.5){0}
  \put(24,5.3){\small $\alpha={1\over 2}$}
 
   \put(7,35){ $\frac 1 d$}
  \put(9, 45){1}
   \put(1,58){\small $\dim_H E_\varphi({\Phi})$}

  \put(80,25){ $1 \over d\beta-\beta+1$}
    \put(20,-1){exponential $e^{n^\alpha}$}
      \put(67,-1){super-exponential $e^{\beta^n}$}
\end{overpic}
\vskip 5pt
\caption{$\dim_H E_{\varphi}(\Phi)$ for $\varphi(j) = {j^a}$.}
\end{figure}

\begin{figure}
\vskip 28pt
\begin{overpic}[scale=1,
]{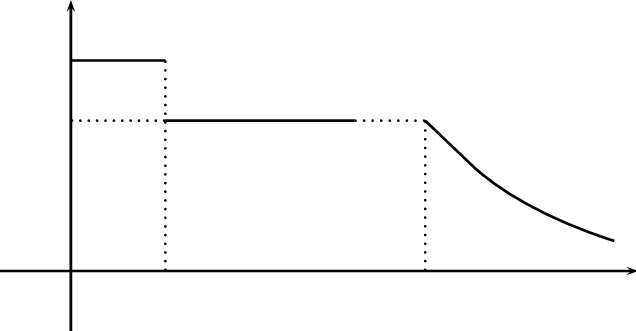}
   \put(64,5.3){$\beta=1+$}
 
  \put(8.8,5.8){0}
  \put(22,5.6){\small $\alpha={b \over b +1}$}
 
  \put(7,32){ $\frac 1 d$}
  \put(9, 42){1}
   \put(1,55){\small $\dim_H E_\varphi({\Phi})$}

  \put(80,25){${1\over d\beta^{1 \over b} -\beta^{1 \over b} +1}$}
    \put(20,-1){exponential $e^{n^\alpha}$}
      \put(67,-1){super-exponential $e^{\beta^n}$}
\end{overpic}
\vskip 5pt
\caption{$\dim_H E_\varphi({\Phi})$ for $\varphi(j) = e^{(\log j)^b}$.}
\end{figure}

\vskip 20pt

\begin{figure}
\vskip 25pt
\begin{overpic}[scale=.6,
]{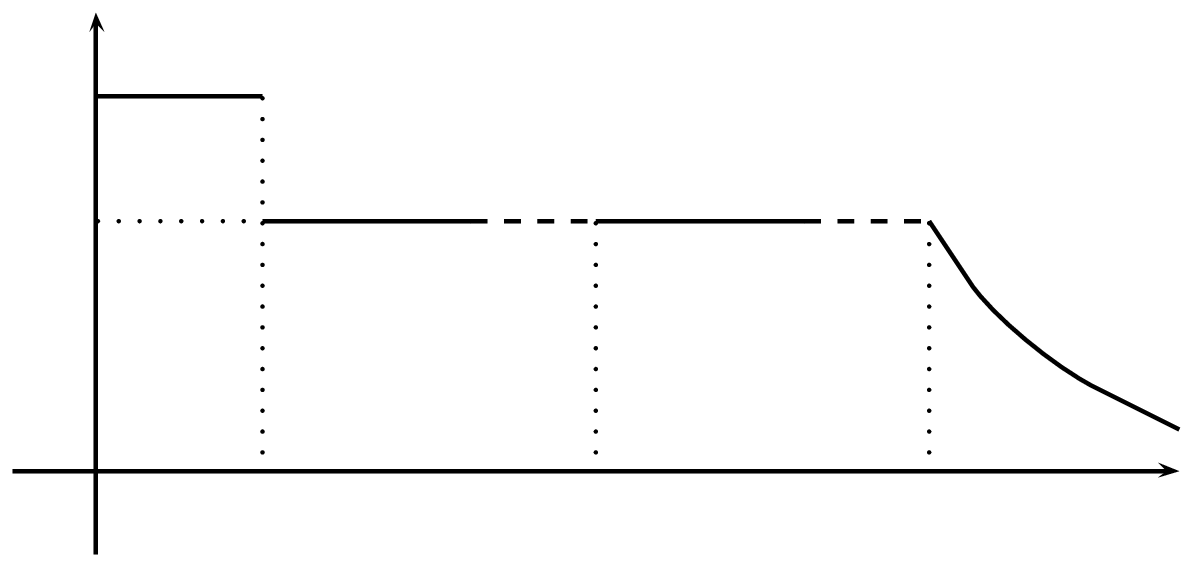}

  \put(47.5,5){$\beta=1+$}
  \put(6,5.5){0}
  \put(19,5.3){\small $\alpha={1\over 2}$}
    \put(74,5.5){\small $\gamma={1+}$}

 
  \put(1,28){ ${1-c \over d}$}
  \put(5.8, 38.5){1}
   \put(1,48){\small $\dim_H E_\varphi({\Phi})$}
   
  \put(82,25){ ${1-c \over d\gamma -(1-c)(\gamma -1)}$}
    \put(16,0){exponential $e^{n^\alpha}$}
      \put(55,0){super-exp $e^{\beta^n}$}
            \put(79,0){sup-sup-exp $e^{e^{\gamma^n}}$}

\end{overpic}
\vskip 5pt
\caption{$\dim_H E_\varphi({\Phi})$ for $\varphi = e^{j^c}$ with $0<c<1$.}
\end{figure}

\begin{figure}
\vskip 25pt
\begin{overpic}[scale=.6,
]{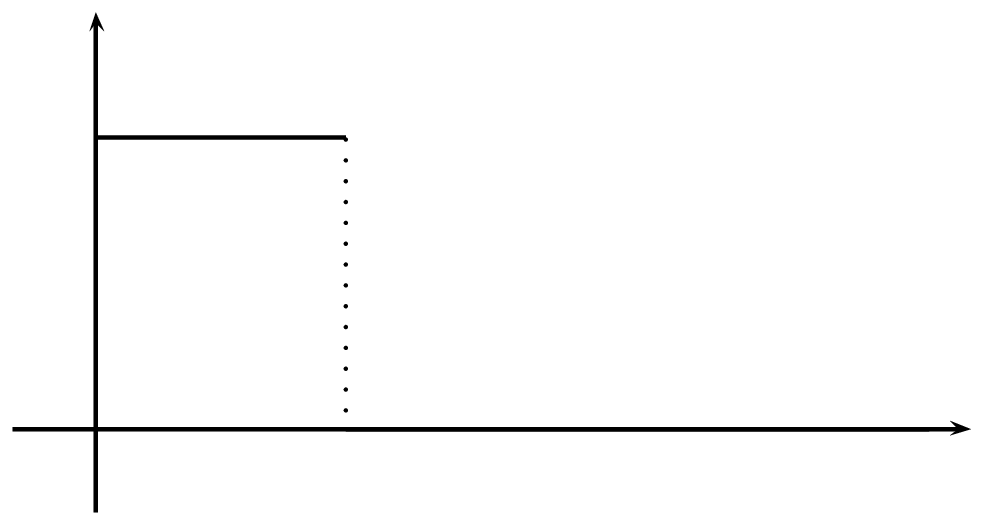}

  \put(30.5,5){$\alpha=1$}
 
  \put(6,5.5){0}
 
  \put(5.8, 38.5){1}
   \put(1,55){\small $\dim_H E_\varphi({\Phi})$}
   
    \put(25,-0.7){exponential $e^{n^\alpha}$}

\end{overpic}
\vskip 5pt
\caption{$\dim_H E_\varphi({\Phi})$ for $\varphi = e^{j^c}$ with $c\geq 1$.}
\end{figure}

\begin{remark}
The critical cases $\alpha={1 \over 2}$ in Theorems \ref{thm-1}, $\alpha={b \over b+1}$ in Theorem \ref{thm-2}, and $\alpha=1$ in Theorems \ref{thm-3} and \ref{thm-4} are not investigated in this paper. However, Theorem 1.2 in \cite{LR16a} suggests that the Hausdorff dimension function has jumps at these points.
\end{remark}

\begin{remark}
Theorem \ref{thm-1} was announced in \cite[Theorem 4.1]{LR16a}, but with an erroneous formula in the part (iii) (now part II).
\end{remark}

\begin{remark}\label{Remark:delta}
For simplicity, in our proofs, we assume $\delta=0$ in the condition (5) of Definition \ref{Def:d-decaying} of the $d$-decaying Gauss like iterated function system. For the general case, the proofs are the same. We need only to replace $d$ by $d+\delta$ for the lower bound and by $d-\delta$ for the upper bound, then take the limit $\delta \to 0$.
\end{remark}

To prove our main theorems, we give four technical lemmas in Section \ref{Sec:Lemmas}.  Lemma \ref{lem-dim-1} is useful to prove full Hausdorff dimension results. Lemma \ref{lem:A} is mainly devoted to proving a lower bound of Hausdorff dimension (sometimes, we can also use it for upper bound). Lemmas \ref{estimation-upper} and \ref{estimation-upper-2} are two combinatorial lemmas serving for the upper bounds of Hausdorff dimension. We believe that our lemmas have independent interests for the future study on multifractal anlysis and Diophantine approximation in infinite iterated function systems or interval maps with infinitely many branches which are often associated to some expansions of numbers. We also stress that though some preliminary versions of our lemmas have already appeared in \cite{WX11}, \cite{Xu}, \cite{FLWW} and \cite{LR16a}, some more efforts are needed to make them applicable for more general settings. Our lemmas are non-trivial generalizations of the corresponding results in \cite{WX11}, \cite{Xu}, \cite{FLWW} and \cite{LR16a}. 

\smallskip
The rest of the paper is organized as follows. In Section \ref{Sec:3}, we give the proofs for (I-1) of Theorems \ref{thm-1}-\ref{thm-4} and (I-2) of Theorem \ref{thm-4}. The remaining proofs are given in the last section. 

\smallskip
{ In the whole paper, for simplicity, for two real sequences $\{f(n)\}_{n\geq 1}$ and $\{g(n)\}_{n\geq 1}$, we will use the symbol $f(n)\approx g(n)$ to denote $f(n)/g(n) = 1+o(1)$, or more precisely
\[
\lim_{n\to\infty}{f(n) \over g(n)}=1.
\]
The symbol $f(n) \asymp g(n)$ means that there exist two positive constants $C_1$ and $C_2$ such that 
\[
C_1\leq {f(n) \over g(n)} \leq C_2.
\]}

\bigskip

\section{Technical lemmas}\label{Sec:Lemmas}

In this section, we will prove four technical lemmas. 
{Before giving the four technical lemmas and their proofs, let us first state a variant of Mass Distribution Principle (\cite[Principle 4.2]{Fa}), sometimes called Billingsley Lemma in the literature, which is a main tool for the lower bound estimation of Hausdorff dimension. 

Denote by $B_r(x)$ the ball with center $x\in \mathbb{R}$ and radius $r>0$.
\begin{lemma}[see \cite{Fa97} Proposition 2.3]\label{lem:MDP}
Let $E\subset \mathbb{R}$ be a Borel set and let $\mu$ be a finite Borel measure.  If the local dimension
\[
\liminf_{r\to 0} {\log \mu(B_r(x)) \over \log r} \geq s
\]
for all $x\in E$ and $\mu(E)>0$, then $\dim_H E \geq s$.
\end{lemma}
}

\medskip

Now, let us go to the four technical lemmas. The first technical lemma serves for the proof of full dimension in the theorems, i.e., the proofs for (I-1) of Theorems \ref{thm-1}-\ref{thm-4}. 

Recalling Remark \ref{Remark:delta}, we remind that we assume $\delta=0$ in the condition (5) of Definition \ref{Def:d-decaying} through all of our proofs.

Let $(n_k)_{k\geq 1}$ be a positive sequence satisfying $n_{k}/k \to \infty$ and $n_{k+1}/n_k \to 1$ as $k\to \infty$. Let $(u_k)_{k\geq 1}$ be a positive sequence such that {$u_k\to \infty$ as $k\to\infty$ and} 
\begin{equation} \label{eqn:usef}
\lim_{k\to\infty} {1 \over n_k} {\sum_{j=1}^k\log u_j }=0.
\end{equation}
For each $M\in \mathbb{N}$, set
\[
E_M:=\{x\in (0,1):  a_{n_k}(x)=u_{k}, \ \text{and } 1\leq a_j(x) \leq M \text{ if } j\neq n_k   \}.
\]
Then we have the following lemma whose idea comes from the proof of \cite[Theorem 1.4]{WX11}.
\begin{lemma}\label{lem-dim-1}
Suppose the $d$-decaying Gauss like iterated function system $\{f_n\}_{n\in\mathbb{N}}$ satisfies the distortion property (\ref{distortion}). 
{Let $(n_k)_{k\geq 1}$, $(u_k)_{k \geq 1}$ and $E_M$ be defined as above.} Then we have 
\[
\lim_{M\to \infty}\dim_H E_M =1.
\]
\end{lemma}
\begin{proof}
For any $k\geq 1$, let $I_{n_k}(a_1\cdots a_{n_k})$ be an $n_k$-basic interval intersecting $E_M$. By the condition (5) of Definition \ref{Def:d-decaying} and the distortion property (\ref{distortion}), we have
\[
|I_{n_k}|\geq K_3^{2k-1}K_1^k \prod_{j=1}^k |I_{n_j-n_{j-1}-1}(a_{n_{j-1}+1}, \cdots, a_{n_j-1})|\cdot a_{n_j}^{-d},
\]
where by convention $n_0=0$.

Let $s(M)$ be the Hausdorff dimension of the set of points $x$ such that all $a_j(x)\leq M$. Then $s(M)$ is increasing to $1$ as $M\to\infty$, see for example, \cite[Theorem 3.15]{MU}. Further, there exists a probability measure $\nu$ living on $[0,1]$ and a positive constant $C_M$ such that for any basic interval $I_n(a_1, \dots, a_n)$ we have
\[
\nu(I_n(a_1, \dots, a_n)) \leq C_M |I_n(a_1, \dots, a_n)|^{2s(M)-1}.
\]

Define a probability measure $\mu$ on each basic interval $I_{n_k}$ intersecting $E_M$ by 
\[
\mu(I_{n_k})=\prod_{j=1}^k \nu(I_{n_j-n_{j-1}-1}(a_{n_{j-1}+1}, \cdots, a_{n_j-1})).
\]
By Kolmogorov Consistence Theorem, $\mu$ is well defined and is supported on $E_M$.

Then for each  $x\in E_M$, we have 
\begin{align}\label{eqn:length_lower}
|I_{n_k}(x)|^{2s(M)-1} \geq (K_3^{2k-1}K_1^k)^{2s(M)-1} C_M^{-k} \mu(I_{n_k}(x)) \left(\prod_{j=1}^k a_{n_j}^{-d}\right)^{{2} s(M)-1}.
\end{align}

Observe that \eqref{eqn:usef} means that $\lim_{k\to\infty}{1 \over n_k}\sum_{j=1}^k \log a_{n_j}=0$, while the part (4) of of Definition \ref{Def:d-decaying} 
implies that 
\begin{equation} \label{eqn:usef2}
\log |I_{n_k}(x)| \leq  \frac {\log A} m n_k.
\end{equation} 
Thus, using \eqref{eqn:length_lower} and \eqref{eqn:usef2}, by simple calculations, we have
\begin{equation} \label{eqn:usef3}
\frac {\log \mu(I_{n_k}(x))} {\log |I_{n_k}(x)|} \geq 2s(M)-1-o(1)
\end{equation}
for large $k$.

This allows us to estimate the local dimension $\liminf_{r \to 0}{\log \mu(B_r(x)) \over \log r}$ of the measure $\mu$ at $x$. 
Let us first observe the following two facts.

\smallskip
{\bf Fact 1}.  Let $x\in E_M$. For $r=|I_{n_k}(x)|$,
\[
B_r(x) \cap E_M \subset I_{n_k}(x).
\]

Indeed, {note that for $x\in E_M$, $a_{n_k}(x)=u_k$ and $f_{u_k}([0,1])=I_1(u_k)$.  By Mean Value Theorem and the condition (5) of Definition \ref{Def:d-decaying}, the basic interval $I_1(u_k)$ has length ${\asymp} u_k^{-d}$ (remind that we assume $\delta=0$).
Further by conditions (2) and (3) of Definition \ref{Def:d-decaying}, on the left of $I_1(u_k)$, there are all $I_1(j)$ with $1\leq j \leq u_k-1$. Thus the distance of $I_1(u_k)$ to $0$ is 
\[
\asymp \sum_{j=1}^{u_k-1} j^{-d} \asymp  1.
\]
While on the right of $I_1(u_k)$, there are all $I_1(j)$ with $j \geq u_k+1$. Thus the distance of $I_1(u_k)$ to $1$ is 
\[
\asymp \sum_{j=u_k+1}^\infty j^{-d} \asymp  u_k^{-d+1}.
\]
 Therefore $I_1(u_k)$ lies in distance larger than ${\asymp} u_k^{-d+1}$ from the endpoints $\{0,1\}$.
 
Now, notice that the pair $(I_{n_k}(x), I_{n_k-1}(x))$ is the image of the pair $(I_1(u_k),[0,1])$ under the map $f_{a_1}\circ\ldots\circ f_{a_{n_k-1}}$. Moreover, by the bounded distortion property (\ref{distortion}), we have 
\[
\sum_{j=1}^{u_k-1}|I_{n_k}[a_1, \cdots ,  a_{n_k-1}, j]|  \asymp \left(\sum_{j=1}^{u_k-1}|I_1(j)| \right) \cdot I_{n_k-1}[a_1 , \cdots ,  a_{n_k-1}],
\]
and
\[
\sum_{u_k+1}^\infty |I_{n_k}[a_1, \cdots , a_{n_k-1}, j]|  \asymp \left(\sum_{u_k+1}^\infty |I_1(j)| \right) \cdot I_{n_k-1}[a_1 , \cdots , a_{n_k-1}].
\]
Thus, the distance of $I_{n_k}(x)$ from the endpoints of $I_{n_k-1}(x)$ is at least a constant multiplied by 
\[
{u_k^{-d+1} \over  u_k^{-d} } |I_{n_k}(x)| = u_k |I_{n_k}(x)|,
\]
which, by the assumption $\lim_{k\to\infty}u_k=\infty$, is much larger than  $3\cdot |I_{n_k}(x)|$. Hence, the ball $B_r(x)$ with $r=|I_{n_k}(x)|$ will never intersect another $(n_k-1)$-basic interval.
}

By the definition of $E_M$, inside $I_{n_k-1}(x)$ there is only one $n_k$-basic interval (which is $I_{n_k}(x)$) intersecting $E_M$. Thus we obtain the assertion of the fact.

\smallskip
{\bf Fact 2}. When $k\to\infty$,
\[
\frac {\log |I_{n_{k+1}}(x)|} {\log |I_{n_{k}}(x)|} \to 1.
\]

Indeed, by the condition (5) of Definition \ref{Def:d-decaying}, {and the distortion property (\ref{distortion}), we have}
\begin{align}\label{eqn:length-quotient}
\frac {|I_{n_{k+1}}(x)|} {|I_{n_{k}}(x)|} \geq {K_3}^{n_{k+1}-n_k} \cdot (K_1M^{-d})^{n_{k+1}-n_k} \cdot K_1 \cdot u_{k+1}^{-d}.
\end{align}
By the hypothesis $n_{k+1}/n_k \to 1$, we have $(n_{k+1}-n_k)/n_k\to 0$.
Then the statement follows from the formulae \eqref{eqn:usef}, \eqref{eqn:usef2} and \eqref{eqn:length-quotient}.

\smallskip

The first fact implies that when $r=|I_{n_k}(x)|$ we can use \eqref{eqn:usef3} in the local dimension calculation. The second fact implies that we do not need to check any $r$ not of the form $r=|I_{n_k}(x)|$. Thus,
by Lemma \ref{lem:MDP}, we have
\[
\dim_H E_M \geq 2s(M)-1.
\]
Passing with $M$ to infinity, we complete the proof of the lemma. 

\end{proof}


\bigskip

The second technical lemma is an improved version of \cite[Lemma 3.2]{FLWW}, \cite[Proof of Theorem 1.3]{JR}, \cite[Lemma 2.2]{LR16a} and  \cite[Lemma 2.2]{LR16b}. 
We often apply it for the lower bounds of Hausdorff dimension. However, we will see that it is also useful for upper bounds. 

Let $\{s_n\}_{n\geq 1}, \{t_n\}_{n\geq 1}$ be two positive real sequences. 
Assume that $s_n>t_n$, $s_n, t_n \to \infty $ as $n\to \infty$, and 
\begin{align}\label{cond:sntn}
\liminf_{n\to \infty}{s_n-t_n \over s_n}>0.
\end{align}
For $N\in \mathbb{N}$, let
\[
B(\{s_n\}, \{t_n\}, N) := \big\{x\in (0,1) : s_n-t_n \leq  {a_n(x)}  \leq s_n+t_n, \ \forall n \geq N\big\}.
\]
We remark that when $n$ is large, by the assumption $t_n\to\infty$, we can always have integers choices for $a_n$ in the interval $[s_n-t_n , s_n+t_n]$. However, for the first terms, there might be no positive integer in the interval $[s_n-t_n , s_n+t_n]$, and hence the set $B(\{s_n\}, \{t_n\}, N) $ is empty.  There is no interest to study such an empty set. Without loss of generality, by modifying the values of finitely many first terms, we always assume that our sequences $\{s_n\}, \{t_n\})$ are such that  $B(\{s_n\}, \{t_n\}, N)$ is nonempty. Further, we will see from our formula of the Hausdorff dimension of $B(\{s_n\}, \{t_n\}, N) $ that the modification of the first finite number of terms of the sequences will not change the Hausdorff dimension.

\begin{lemma} \label{lem:A}
We have 
\[
\dim_H B(\{s_n\}, \{t_n\}, N)  = \liminf_{\ell\to \infty}{\sum_{i=1}^\ell \log t_i \over d \sum_{i=1}^{\ell+1} \log s_i -\log t_{\ell+1}}.
\]
\end{lemma}
\begin{proof}
Within this proof, we write {$f(n)\stackrel{{e}}{\sim}g(n)$} if $f(n)$ and $g(n)$ differ by at most an exponential factor, that is
\[
\limsup_{n\to\infty} {1\over n}\left|\log {f(n) \over g(n)} \right|<\infty.
\]
We give the proof for the case $N=1$. For the general case, note that  
\[
B(\{s_n\}, \{t_n\}, N)=\bigcup_{a_1\cdots a_{N-1}\in \mathbb{N}^{N-1}} f_{a_1}\circ \cdots \circ f_{a_{N-1}} (B(\{s_{n+N-1}\}, \{t_{n+N-1}\}, 1))
\]
 is a countable union of bi-Lipschitz images of $B(\{s_{n+N-1}\}, \{t_{n+N-1}\}, 1)$. Since the bi-Lipschitz maps preserve the Hausdorff dimension, we have 
 \[
 \dim_HB(\{s_n\}, \{t_n\}, N)=\dim_H B(\{s_{n+N-1}\}, \{t_{n+N-1}\}, 1).
 \]
 On the other hand, notice that the dimensional formula of the lemma we will obtain does not depend on the finite number of first terms of the two sequences $\{s_n\}$ and $\{t_n\}$, we then have  
  \[
 \dim_HB(\{s_n\}, \{t_n\}, N)=\dim_H B(\{s_{n}\}, \{t_{n}\}, 1).
 \]

Fix $\ell\geq 1$. Let $I_\ell(a_1, \dots, a_\ell)$ be an $\ell$-basic interval with nonempty intersection with $B(\{s_n\}, \{t_n\}, 1)$. Then for each $1\leq k \leq \ell$, $a_k \in [s_k-t_k, s_k+t_k]$.
 Define 
\[
D_\ell(a_1, \dots, a_\ell):=\big\{x\in I_\ell(a_1, \dots, a_\ell):  a_{\ell+1}(x)\in [s_{\ell+1}-t_{\ell+1}, s_{\ell+1}+t_{\ell+1}] \big\}.
\]
We have 
\begin{align*}
B(\{s_n\}, \{t_n\},1) =& \bigcap_{\ell=1}^\infty \bigcup_{\substack{ a_1,\ldots,a_\ell \\ a_i \in [s_i-t_i, s_i+t_i] }} I_\ell(a_1,\ldots,a_\ell) \\
= &\bigcap_{\ell=1}^\infty \bigcup_{\substack{ a_1,\ldots,a_\ell \\ a_i \in [s_i-t_i, s_i+t_i] }} D_\ell(a_1,\ldots,a_\ell).
\end{align*}
At level $\ell$, we have $\stackrel{{e}}{\sim} \prod_{i=1}^\ell t_i$ intervals $I_\ell(a_1,\ldots,a_\ell)$ and thus $\stackrel{{e}}{\sim} \prod_{i=1}^\ell t_i$ corresponding $D_\ell(a_1,\ldots,a_\ell)$. By the condition (5) of Definition \ref{Def:d-decaying}, and the assumption \eqref{cond:sntn}, each $I_\ell(a_1,\ldots,a_\ell)$ is of size $\stackrel{{e}}{\sim} \prod_{i=1}^\ell s_i^{-d}$. Moreover,
\[
\frac {|D_\ell(a_1,\ldots,a_\ell)|}{|I_\ell(a_1,\ldots,a_\ell)|} \ \stackrel{{e}}{\sim} \ \sum_{i=s_{\ell+1}-t_{\ell+1}}^{s_{\ell+1}+t_{\ell+1}} i^{-d}  \ \stackrel{{e}}{\sim} \ t_{\ell+1} s_{\ell+1}^{-d}.
\]
Thus, using for a given $\ell$ the sets $D_\ell(a_1,\ldots,a_\ell)$ as a cover for $B(\{s_n\}, \{t_n\}, 1)$, we need $\stackrel{{e}}{\sim} \prod_{i=1}^\ell t_i$ of them, each of size $\stackrel{{e}}{\sim} t_{\ell+1} \prod_{i=1}^{\ell+1} s_i^{-d}$.  Therefore, by regular calculation, we can obtain the upper bound
\[
\dim_H B(\{s_n\}, \{t_n\}, 1) \leq \liminf_{\ell\to\infty} \frac {\sum_{i=1}^\ell \log t_i}{d\sum_{i=1}^{\ell+1} \log s_i - \log t_{\ell+1}}.
\]

\medskip

To get the lower bound, we consider a probability measure $\mu$ uniformly distributed on $ B(\{s_n\}, \{t_n\}, 1)$, in the following sense. 
{By the definition of $B(\{s_n\}, \{t_n\}, 1)$, when $a_1, \ldots, a_{\ell-1}$ are fixed, $a_\ell$ can only take finitely many values $s_\ell-t_\ell,\ldots, s_\ell+t_\ell$. Then the conditional probability of every possible value of $a_\ell$ (conditioned on $a_1,\ldots, a_{\ell-1}$) is chosen to be equal, that is equal to $1/(2t_\ell+1)$.}
The basic intervals $I_\ell(a_1,\ldots,a_\ell)$ and corresponding $D_\ell(a_1,\ldots,a_\ell)$ have the measure $\stackrel{{e}}{\sim} \prod_{i=1}^\ell t_i^{-1}$. 

Our goal is to apply Lemma \ref{lem:MDP}, hence we need to calculate the local dimension of the measure $\mu$ at a $\mu$-typical point $x\in B(\{s_n\}, \{t_n\},1)$. 
Fix any $x\in B(\{s_n\}, \{t_n\},1)$. Denote by $r_\ell$ the diameter of the set $D_\ell(a_1(x),\ldots,a_\ell(x))$ and by $r_\ell'$ the diameter of $I_\ell(a_1(x),\ldots,a_\ell(x))$. When $r=r_\ell$, we have
\[
\frac {\log \mu(B_r(x))} {\log r} =\frac {\log \mu(B_{r_\ell}(x))} {\log r_\ell} =\frac {\log \mu(D_\ell(a_1(x),\ldots,a_\ell(x)))} {\log r_\ell}.
\]
Since $\mu(D_\ell(a_1,\ldots,a_\ell))\stackrel{{e}}{\sim} \prod_{i=1}^\ell t_i^{-1}$ and $r_\ell\stackrel{{e}}{\sim} t_{\ell+1} \prod_{i=1}^{\ell+1} s_i^{-d}$, we have
\begin{align}\label{eqn:local-dimension-rn}
\liminf_{\ell\to\infty}\frac {\log \mu(B_{r_\ell}(x))} {\log r_\ell}  =\liminf_{\ell\to\infty} \frac {\sum_{i=1}^\ell \log t_i} {d\sum_{i=1}^{\ell+1} \log s_i - \log t_{\ell+1}}.
\end{align}

For $r_\ell<r<r_\ell'$, the ball $B_r(x)$ still does not intersect any point from $B(\{s_n\}, \{t_n\},1)\setminus D_\ell(a_1(x),\ldots,a_\ell(x))$, hence it has the same measure as $B_{r_\ell}(x)$, but a larger diameter.  Thus, for $r_\ell<r<r_\ell'$,
\[
\frac {\log \mu(B_{r_\ell}(x))} {\log r_\ell} < \frac{\log \mu(B_{r}(x))} {\log r}.
\]

Finally, since each basic interval $I_{\ell+1}(a_1(x), \dots, a_\ell(x),j)$ contained in $D_\ell(a_1(x), \dots, a_\ell(x))$ has  the same measure and approximately the same diameter, we have for $r_{\ell+1}'<r<r_\ell$,
{\[
\mu(B_r(x)) \leq {2 r  \over r_{\ell+1}' }  \cdot {K_2 \over K_1}({s_{\ell+1}+ t_{\ell+1} \over s_{\ell+1}- t_{\ell+1}})^d \cdot {\mu (B_{r_\ell}(x)) \over 2 t_{\ell+1}}.
\]
Since  
\[
 {2r_\ell \over r'_{\ell+1} 2 t_{\ell+1}}\leq 2 \cdot {K_2 \over K_1}({s_{\ell+1}+ t_{\ell+1} \over s_{\ell+1}- t_{\ell+1}})^d,
\]
we have 
\[
\mu(B_r(x)) \leq 2  \cdot \left({K_2 \over K_1}({s_{\ell+1}+ t_{\ell+1} \over s_{\ell+1}- t_{\ell+1}})^d\right)^2 \cdot {\mu (B_{r_\ell}(x))\cdot {r \over r_\ell} }.
\]
By (\ref{cond:sntn}), there is a constant $K_6>0$, such that 
\[
2  \cdot \left({K_2 \over K_1}({s_{\ell+1}+ t_{\ell+1} \over s_{\ell+1}- t_{\ell+1}})^d\right)^2 \leq K_6.
\]
Applying the obvious fact that
\[
\frac {\log z_1 z_2}{\log z_1 z_3} > \frac {\log z_2}{\log z_3}
\]
for all $z_1<1$ and $z_3<z_2<1$, we see that for $r_{\ell+1}'<r<r_\ell$,
\[
\frac {\log \mu(B_{r_\ell}(x))} {\log r_\ell} < \frac{\log (\mu(B_{r_\ell}(x)) \cdot r/r_\ell)} {\log r} \leq  \frac{\log \mu(B_{r}(x))- \log K_6} {\log r}.
\]}

Thus, the minimum of the function $r\to \log \mu(B_r(x))/\log r$ for $r_{\ell+1}'<r<r_\ell'$ is equal to its value at $r_\ell$, up to an error term that vanishes as $\ell\to\infty$. 
Therefore, 
\[
\liminf_{r\to 0} \frac {\log \mu(B_r(x))} {\log r} =\liminf_{\ell\to \infty} \frac {\log \mu(B_{r_\ell}(x))} {\log r_\ell}. 
\]
Applying Lemma \ref{lem:MDP}, by \eqref{eqn:local-dimension-rn}, we obtain the lower bound 
\[
\dim_H B(\{s_n\}, \{t_n\}, 1) \geq \liminf_{\ell\to\infty} \frac {\sum_{i=1}^\ell \log t_i}{d\sum_{i=1}^{\ell+1} \log s_i - \log t_{\ell+1}}.
\]
The proof is thus completed.
\end{proof}

\medskip
The remaining two technical lemmas are generalizations of \cite[Lemma 2.1]{LR16a}. They are the key to prove the upper bounds of Hausdorff dimension.
In the proofs, the Riemann zeta function $\zeta(\cdot)$ will be often used. 

For $m,n\in \mathbb{N}$, $a>0$ and $\epsilon>0$, 
let 
\[
A(m,n, a,\epsilon):=\left\{(i_1, \dots, i_n)\in \mathbb{N}^n: \ \sum_{k=1}^n i_k^a\in [m, m(1+\epsilon)]\right\}.
\]
For $s>1/d$, write
 \[
 G(m,n,a,\epsilon, s)= \sum_{i_1\cdots i_n\in A(m,n,a,\epsilon)} \prod_{k=1}^ni_k^{-ds}.
 \]

\begin{lemma}\label{estimation-upper}
There exist positive constants $C_1=C_1(a,s), C_2=C_2(s)$, and $C_3=C_3(a)$, such that for all $C_3 \cdot  (m3^{2-n})^{-1/a}<\varepsilon\leq1/3$, we have
\[
  G(m,n, a, \epsilon, s) \leq C_1C_2^{n-1}\epsilon m^{1-ds \over a}.
\]
\end{lemma}
\begin{proof}

The proof goes by induction. First consider the case $n=2$. 
Note that if $i_1^a+i_2^a \in [m, m(1+\epsilon)] $ then at most one of $i_1^a,  i_2^a$ is strictly larger than ${m(1+\epsilon) \over 2}$. We divide the sum in the definition of $G(m,n,a,\epsilon, s)$ into two parts, one is $i_1^a\leq {m(1+\epsilon)\over 2}$, the other is $i_1^a> {m(1+\epsilon) \over 2}$. However, by permuting $i_1$ and $i_2$, the latter is smaller than the former. Thus, 
\begin{align*}
   G(m,2, a, \epsilon,s)
    \leq 2\sum_{k=1}^{({m(1+\varepsilon) \over 2})^{{1 \over a}}} k^{-ds}(m-k^a)^{-{ds\over a}}  \cdot N_{m,a,\epsilon}(k),
\end{align*}
with $N_{m,a,\epsilon}(k):=\sharp \{i_2: m-k^a \leq i_2^a \leq m(1+\epsilon)-k^a\}$.

Assuming $\varepsilon\leq 1/3$, we can estimate for $a\geq 1$
\begin{align*}
 N_{m,a,\epsilon}(k)
 \leq \lceil a^{-1}\varepsilon m (m-k^a)^{\frac 1a -1}\rceil 
 \leq  \lceil \varepsilon m^{1/a} \cdot a^{-1}  3^{1-1/a} \rceil,
\end{align*}
while for $a<1$ 
\begin{align*}
 N_{m,a,\epsilon}(k)
\leq  \lceil a^{-1}\varepsilon m (m(1+\varepsilon)-k^a)^{\frac 1a -1}\rceil 
\leq  \lceil \varepsilon m^{1/a} \cdot a^{-1}  (4/3)^{1/a-1} \rceil.
\end{align*}
That is, in both cases we will get an upper estimation in the form $$\lceil  \varepsilon m^{1/a} \cdot C_4(a) \rceil.$$

If $z>1$, we can estimate $\lceil z \rceil \leq 2z$. Thus, for $\varepsilon>1/(m^{1/a} \cdot C_4(a))$ we have 
\[
 N_{m,a,\epsilon}(k) \leq 2\varepsilon m^{1/a} \cdot C_4(a).
\]
Hence,
\begin{align}\label{formula-G-2}
\begin{split}
& G(m,2, a, \epsilon,s)\\
 \leq & 2 \sum_{k=1}^{({m(1+\epsilon) \over 2})^{{1 \over a}}} k^{-ds} ({m\over 3})^{-{ds\over a}} \cdot 2\varepsilon  m^{{1\over a}}\cdot C_4(a) \\
 \leq &4 \cdot \zeta(ds) \cdot  3^{\frac {ds} a } \cdot C_4(a)  \cdot \varepsilon m^{\frac {1-ds} a}.
 \end{split}
\end{align}
Let
\[
C_1=  2\cdot 3^{\frac {ds} a} C_4(a), \  C_2= 6 \cdot 3^{ds-1 \over a} \cdot \zeta(ds).
\]
Then by (\ref{formula-G-2}), we have
\[
G(m,2, a, \epsilon,s)  \leq C_1C_2 \epsilon   m^{1-ds \over a}.
\]

Assume now that the assertion is satisfied for all $n<N$ for some $N>2$, we will prove by induction that it holds for $n=N$ as well. 

As above, there is at most one $i_k$ such that $i_k^a>{m(1+\epsilon) \over 2}$. Thus the sum of $G(m,N,a,\varepsilon, s)$ can be divided into two parts, one is $i_1^a\leq {m(1+\epsilon) \over 2}$ 
and the other is $i_1^a>{m(1+\epsilon) \over 2}$. Again, the latter is smaller than the former because we can permute $i_1$ and $i_2$. Thus
\begin{align*}
G(m,N,a,\varepsilon, s) \leq 2 \sum_{k=1}^{({m(1+\epsilon) \over 2})^{{1 \over a}}} k^{-ds}  \sum_{i_2\cdots i_n \in \widetilde{A}(m,n, k, a,\epsilon)} \prod_{k=2}^ni_k^{-ds},
\end{align*}
where
\[
\widetilde{A}(m,n, k, a,\epsilon):=\left\{(i_2, \dots, i_n)\in \mathbb{N}^n: \ \sum_{k=2}^n i_k^a\in [m-k^a, m(1+\epsilon)-k^a]\right\}.
\]

Further, observe  that $3(m-k^a)\epsilon \geq m\epsilon$. Then the sum range $$[m-k^a, \ m(1+\epsilon)-k^a]$$ in $\widetilde{A}(m,n, k, a,\epsilon)$ is covered by the union of 
$$[m-k^a, \ (m-k^a)(1+\epsilon)], \ [(m-k^a)(1+\epsilon), \ (m-k^a)(1+\epsilon)^2], $$
and
$$[(m-k^a)(1+\epsilon)^2, \ (m-k^a)(1+\epsilon)^3].$$
Hence,
\begin{align*}
G(m,N,a,\varepsilon, s) \leq 2 \sum_{k=1}^{({m(1+\epsilon) \over 2})^{{1 \over a}}} k^{-ds} \sum_{j=0}^2 G\big((m-k^a)(1+\epsilon)^j, N-1, a, \varepsilon, s\big).
\end{align*}
Note that
\begin{align}\label{estimate: m-ka}
(m-k^a)(1+\epsilon)^j \geq {m \over 3}, \ \text{for} \ j=0,1,2.
\end{align}
Substituting the induction assumption, we get
\begin{align*}
G(m,N,a,\varepsilon, s)& \leq 6 \cdot C_1C_2^{N-2} \varepsilon ({m \over 3})^{1-ds\over a} \sum_{k=1}^{({m(1+\epsilon) \over 2})^{{1 \over a}}} k^{-ds} \\
&\leq 6 \cdot 3^{ds-1 \over a} C_1C_2^{N-2} \varepsilon m^{1-ds\over a} \zeta(ds)\\
&=C_1C_2^{N-1} \varepsilon m^{1-ds\over a},
\end{align*}
where the last equality comes from the definition of $C_2$.
%
Remember, however, when we prove our lemma for the case $n=2$, we have assumed 
\[
{1 \over m^{1/a} C_4(a)} < \epsilon \leq  1/3.
\] 
So, by \eqref{estimate: m-ka}, when we conduct the induction we need to assume that 
\[
{1 \over ({m \over 3^i})^{1/a} C_4(a)} < \epsilon \leq 1/3, \ \text{for} \ i=0, \dots, N-2.
\] 
That is, we need to assume 
$$\varepsilon \in \big((m 3^{2-N})^{-1/a} (C_4(a))^{-1},\ 1/3\big].$$ Taking  $C_3=1/C_4(a)$, we finish the proof.
\end{proof}

\medskip
The next lemma is very similar. Let
\[
\widehat{A}(m,n, b,\epsilon):=\left\{(i_1, \dots, i_n)\in \mathbb{N}^n: \ \sum_{k=1}^n e^{(\log i_k)^b}\in [m, m(1+\epsilon)]\right\},
\]
and for $s>1/d$, write
\[
\widehat{G}(m,n,b,\epsilon, s)= \sum_{i_1\cdots i_n\in \widehat{A}(m,n,b,\epsilon)} \prod_{k=1}^n i_k^{-ds}.
\]

\begin{lemma}\label{estimation-upper-2}
There exists a positive constant $\widehat{C}=\widehat{C}(s)$ such that for all $e^{-(\log (m3^{2-n}))^{1/b}}<\varepsilon \leq 1/3$, we have
\[
  \widehat{G}(m,n, b, \epsilon, s) \leq 6\cdot \widehat{C}^{n-1}\epsilon \cdot  e^{(1-ds)(\log m)^{1/b}} .
\]
\end{lemma}
\begin{proof}

The proof goes again by induction. First, consider the case $n=2$. 
Similar to the proof of Lemma \ref{estimation-upper}, we have 
\begin{align*}
   \widehat{G}(m,2, b, \epsilon,s)
    \leq 2\sum_{k=1}^{e^{(\log (m(1+\varepsilon)/2)^{1/b}}} k^{-ds}e^{- ds(\log (m-e^{(\log k)^b}))^{1/b}} \cdot \widehat{N}_{m,b,\epsilon}(k), 
\end{align*}
with $$ \widehat{N}_{m,b,\epsilon}(k):=\sharp \Big\{i_2: m-e^{(\log k)^b} \leq e^{(\log i_2)^b}\leq m(1+\epsilon)-e^{(\log k)^b} \Big\}.$$

For $\epsilon \leq   1/3$, short calculations give us the following estimation
\[
\widehat{N}_{m,b,\epsilon}(k) \leq \lceil 3\epsilon \cdot e^{(\log m)^{1/b}}  \rceil.
\]
Hence, if $\epsilon > e^{-(\log m)^{1/b}}$,
\[
\widehat{N}_{m,b,\epsilon}(k) \leq 6 \epsilon \cdot e^{(\log m)^{1/b}} .
\]
Thus, by noting $e^{\log (m/3))^{1/b}} \geq {1\over 3} e^{(\log m)^{1/b}}$, we obtain
\[
 \widehat{G}(m,2, b, \epsilon,s) \leq 12 \cdot 3^{ds}\zeta(ds) e^{(1-ds)(\log m)^{1/b}} .
\]

Assume now that the assertion is satisfied for all $n<N$ for some $N>2$, we will prove by induction that it holds for $n=N$ as well. We have
\[
\widehat{G}(m,N,b,\varepsilon, s) \leq 2 \sum_{k=1}^{e^{(\log (m(1+\varepsilon)/2)^{1/b}}} k^{-ds} \sum_{j=0}^2 \widehat{G}((m-e^{(\log k)^b})(1+\epsilon)^j, N-1, b, \varepsilon, s).
\]
Substituting the induction assumption, we get
\[
\widehat{G}(m,N,b,\varepsilon, s) \leq 12 \cdot 3^{ds} \widehat{C}^{N-2} \varepsilon e^{(1-ds)(\log m)^{1/b}}  \zeta(ds).
\]

Thus, we have proved the assertion for 
\[
\widehat{C}=2 \cdot 3^{ds} \zeta(ds)
\]
under the assumption $$\varepsilon\in \big(e^{-(\log (m3^{2-N}))^{1/b}},\ 1/3 \big].$$

\end{proof}

\bigskip

\section{Proofs for (I-1) of Theorems \ref{thm-1}-\ref{thm-4} and (I-2) of Theorem \ref{thm-4}}\label{Sec:3}

\subsection{Proofs for (I-1) of Theorems \ref{thm-1}-\ref{thm-4}}

For these parts of proofs we suppose that the $d$-decaying Gauss like iterated function system satisfies the distortion property (\ref{distortion}). We will apply Lemma \ref{lem-dim-1}.

Note that in all cases we are going to prove, the function $\Phi$ is taken as $\Phi(n)=e^{n^\alpha}$. Let $\epsilon>0$. Take $n_k=k^{{1 \over \alpha}(1-\epsilon)}$ and $u_k=\varphi^{-1}(\Phi(n_k)-\Phi(n_{k-1}))$. 

The sequence $(n_k)_{k\geq 1}$ clearly satisfies $$\lim_{k\to\infty} n_k/k\to\infty \ \text{and} \ \lim_{k\to\infty} n_{k+1}/n_k=1,$$ which we will need to eventually apply Lemma \ref{lem-dim-1}. 

We can also check that $E_M \subset E_\varphi(\Phi)$. In fact, for any $x\in E_M$ we have
\[
\Phi(n_k) < S_{n_k}\varphi(x) < \Phi(n_k) + n_k \varphi(M).
\]
Since $\Phi(n)/n\to \infty$, we see that 
\[
\frac {S_{n_k}\varphi(x)} {\Phi(n_k)} \to 1.
\]
However, as $n_{k+1}/n_k \to 1$ and $S_n\varphi$ is increasing, this is enough to have
\[
\lim_{n\to\infty} \frac {S_{n}\varphi(x)} {\Phi(n)} = \lim_{k\to\infty} \frac {S_{n_k}\varphi(x)} {\Phi(n_k)}
\]
and we are done.

Now we need only to check for each case of $\varphi$ in Theorems \ref{thm-1}-\ref{thm-4}, {$u_k\to\infty$ as $k\to \infty$,} and the condition (\ref{eqn:usef}) is satisfied.
First, notice that
\[
\Phi(n_k)-\Phi(n_{k-1})= e^{k^{1-\epsilon}}-e^{(k-1)^{1-\epsilon}}. 
\]
Thus, when $\varphi(j)=j^a$, by Mean Value Theorem, we have 
\[
\big((1-\epsilon)k^{-\epsilon}e^{(k-1)^{1-\epsilon}}\big)^{1/a} \leq u_k \leq  \big((1-\epsilon)(k-1)^{-\epsilon}e^{k^{1-\epsilon}}\big)^{1/a},
\]
{and $\lim_{k\to\infty}u_k=\infty$. Further,} if $\alpha<1/2$ and $\epsilon$ is small enough, 
\[
\lim_{k\to\infty}{1 \over n_k} {\sum_{j=1}^k \log u_j}= \lim_{k\to\infty} {\sum_{j=1}^k j^{1-\epsilon}/a \over k^{{1 \over \alpha}(1-\epsilon) }}=0. 
\]
When $\varphi(j)=e^{(\log j)^b}$, then 
\[
e^{(\log ((1-\epsilon)k^{-\epsilon}e^{(k-1)^{1-\epsilon}}))^{1/b}} \leq u_k \leq e^{(\log ((1-\epsilon)(k-1)^{-\epsilon}e^{k^{1-\epsilon}}))^{1/b}},
\]
{and $\lim_{k\to\infty}u_k=\infty$. Further,}  if $\alpha<{b \over b+1}$, and $\epsilon$ is small enough, 
\[
\lim_{k\to\infty}{1 \over n_k} {\sum_{j=1}^k \log u_j}= \lim_{k\to\infty} {\sum_{j=1}^k j^{1-\epsilon \over b} \over k^{{1 \over \alpha}(1-\epsilon) }}=0. 
\]
When $\varphi(j)=e^{j^c}$, we have
\[
 \log ((1-\epsilon)k^{-\epsilon}e^{(k-1)^{1-\epsilon}})^{1/c} \leq u_k \leq \log ((1-\epsilon)(k-1)^{-\epsilon}e^{k^{1-\epsilon}})^{1/c},
\]
{and $\lim_{k\to\infty}u_k=\infty$. Further,} if $\alpha<1$ and $\epsilon$ is small enough, 
\[
\lim_{k\to\infty}{1 \over n_k} {\sum_{j=1}^k \log u_j}= \lim_{k\to\infty} {\sum_{j=1}^k {1-\epsilon \over c}\log j \over k^{{1 \over \alpha}(1-\epsilon) }}=0. 
\]
Hence, in all cases, {$\lim_{k\to\infty}u_k=\infty$} and the condition (\ref{eqn:usef}) is satisfied. 

Applying Lemma \ref{lem-dim-1}, we complete the proofs.

\medskip

\subsection{Proof for (I-2) of Theorem \ref{thm-4}}

We will use a natural covering. Suppose $\Phi(n)=e^{n^\alpha}$ with $\alpha>1$. 
For each $x\in E_\varphi(\Phi)$, for any small $\epsilon>0$, and for all large enough $n$, we have
\[
(1-\epsilon) \Phi(n) \leq \sum_{k=1}^n \varphi(a_k) \leq (1+\epsilon)\Phi(n).
\]
Thus 
\[
(1-\epsilon)\Phi(n)- (1+\epsilon) \Phi(n-1)  \leq \varphi(a_n)\leq (1+\epsilon) \Phi(n)-(1-\epsilon)\Phi(n-1).
\]
Note that for $\alpha>1$, we have 
\[
(1+\epsilon) \Phi(n)-(1-\epsilon)\Phi(n-1)= (1+\epsilon)e^{n^\alpha} - (1-\epsilon)e^{(n-1)^\alpha} \leq  (1+\epsilon)e^{n^\alpha},
\]
and 
\[
(1-\epsilon) \Phi(n)-(1+\epsilon)\Phi(n-1)= (1-\epsilon)e^{n^\alpha} - (1+\epsilon)e^{(n-1)^\alpha} \geq  (1-2\epsilon)e^{n^\alpha}.
\]
Then 
\[
(1-2\epsilon)e^{n^\alpha} \leq \varphi(a_n)\leq (1+\epsilon)e^{n^\alpha}.
\]
However, for $\varphi(j)=e^{j^c}$ with $c\geq 1$, there is at most one $j$ such that
\[
(1-2\epsilon)e^{n^\alpha} \leq \varphi(j) \leq (1+\epsilon)e^{n^\alpha}.
\]
Hence, $E_\varphi(\Phi)$ is a countable set which has Hausdorff dimension $0$.

\bigskip

\section{Remaining proofs}

We will divide the case I-2 of Theorem \ref{thm-1} into two subcases: subcase I-2a for $1/2<\alpha<1$, and subcase I-2b for $\alpha\geq 1$. Similarly, we will divide the case I-2 of Theorem \ref{thm-2} into subcase I-2a ($b/(b+1)<\alpha<1$) and subcase I-2b ($\alpha\geq 1$). 

Theorem \ref{thm-1}, case II; Theorem \ref{thm-1}, subcase I-2b; Theorem \ref{thm-2}, case II; Theorem \ref{thm-2}, subcase I-2b; Theorem \ref{thm-3}, case I-2;  Theorem \ref{thm-3}, case II; Theorem \ref{thm-3}, case III are all obtained by applying  Lemma \ref{lem:A}.

\subsection{ Proof of Theorem \ref{thm-1}, case II}

Let $x\in E_\varphi(\Phi)$. Fix some small $\varepsilon >0$. Then there exists $N\in \mathbb{N}$ such that for all $n>N$, $$\Phi(n)(1-\varepsilon) < S_n \varphi(x) < \Phi(n)(1+\varepsilon).$$ This implies
\begin{align}\label{eqn:diff}
\begin{split}
& \varphi(a_n(x)) = S_n \varphi(x) - S_{n-1}\varphi(x) \\
\in& \Big(\Phi(n)(1-\varepsilon)-\Phi(n-1)(1+\varepsilon), \ \Phi(n)(1+\varepsilon)-\Phi(n-1)(1-\varepsilon)\Big)
\end{split}
\end{align}
for $n\geq N$. Substituting the formula for $\Phi$, we get
\[
\varphi(a_n(x)) \in \Big(e^{\beta^n}(1-2\varepsilon), \ e^{\beta^n}(1+2\varepsilon)\Big).
\]
Hence, a further substitution of the formula for $\varphi$ gives us
\[
e^{\beta^n/a} (1-3\varepsilon/a) < a_n(x) < e^{\beta^n/a}(1+3\varepsilon/a).
\]
Put $s_n=e^{\beta^n/a}$ and $t_n=3\varepsilon e^{\beta^n/a}/a$. Then
\[
E_\varphi(\Phi) \subset \bigcup_N B(\{s_n\}, \{t_n\},N).
\]
  By Lemma \ref{lem:A}, we have the upper bound 
\begin{align*}
 \dim_H E_\varphi(\Phi) \leq &\liminf_{\ell\to \infty}{\sum_{j=1}^\ell \log 3\varepsilon e^{\beta^j/a}/a \over d \sum_{j=1}^{\ell+1} \log e^{\beta^j/a} -\log 3\varepsilon e^{\beta^{\ell+1}/a}/a}\\
 =& \liminf_{\ell\to \infty}{\sum_{j=1}^\ell {\beta^j/a} \over d \sum_{j=1}^{\ell+1} {\beta^j/a} -{\beta^{\ell+1}/a}}\\
 =&{1 \over d\beta-\beta+1}.
\end{align*}

On the other hand, let $\varepsilon_n$ be a sequence of positive numbers converging to $0$. Let $x\in B(e^{\beta^n/a}, \varepsilon_n e^{\beta^n/a},1)$. For large $n$ we have 
\[
e^{\beta^n}(1-\varepsilon_n)^a < S_n \varphi(x) < e^{\beta^n}(1+\varepsilon_n)^a + \sum_{i=1}^{n-1} (1+\varepsilon_i)^a \cdot e^{\beta^i} < e^{\beta^n}(1+a\varepsilon_n+o(1)).
\]
Thus, 
\[
E_\varphi(\Phi) \supset  B(e^{\beta^n/a}, \ \varepsilon_n e^{\beta^n/a}, \ 1).
\]
Applying Lemma \ref{lem:A} and doing almost the same calculation as above, we obtain the lower bound.

\subsection{Theorem \ref{thm-1}, case I-2b}
We can repeat the proof of Theorem \ref{thm-1}, case II. From the formula (\ref{eqn:diff}), we get
\[
\varphi(a_n(x)) \in \Big(e^{n^\alpha}(1-2\varepsilon), \ e^{n^\alpha}(1+2\varepsilon)\Big).
\]
Hence,
\[
E_\varphi(\Phi) \subset \bigcup_N B(e^{n^\alpha/a}, \ 3\varepsilon e^{n^\alpha/a}/a,\ N).
\]
On the other hand, for a sequence of positive numbers $\epsilon_n$ converging to $0$, we have
\[
E_\varphi(\Phi) \supset  B(e^{n^\alpha/a}, \ \varepsilon_n e^{n^\alpha/a}, \ 1).
\]
Applying Lemma \ref{lem:A}, we have
\[
\dim_H E_\varphi(\Phi) = \liminf_{\ell\to\infty} \frac {\sum_{j=1}^\ell j^\alpha/a}{d\sum_{j=1}^{\ell+1} j^\alpha /a- (\ell+1)^\alpha/a} = \frac 1d.
\]

\subsection{Theorem \ref{thm-2}, case II} 
From the formula (\ref{eqn:diff}), we get
\[
\varphi(a_n(x)) \in \Big(e^{\beta^n}(1-2\varepsilon), \ e^{\beta^n}(1+2\varepsilon)\Big).
\]
Hence,
\[
E_\varphi(\Phi) \subset \bigcup_N B\Big(e^{\beta^{n/b}}, \ {3\varepsilon \over b} \beta^{n (1/b-1)} e^{\beta^{n/b}}, \ N\Big).
\]
On the other hand, for a positive sequence $\epsilon_n$ converging to $0$, we have
\[
E_\varphi(\Phi) \supset  B\big(e^{\beta^{n/b}}, \ \varepsilon_n \beta^{n (1/b-1)} e^{\beta^{n/b}}, \ 1\big).
\]
Applying Lemma \ref{lem:A}, we have
\[
\dim_H E_\varphi(\Phi) = \liminf_{\ell\to\infty} \frac {\sum_{j=1}^\ell \beta^{j/b}}{d\sum_{j=1}^{\ell+1} \beta^{j/b} - \beta^{(\ell+1)/b}} 
= \frac 1{d\beta^{1/b}-\beta^{1/b}+1}.
\]

\subsection{Theorem \ref{thm-2}, case I-2b}
From the formula (\ref{eqn:diff}), we get
\[
\varphi(a_n(x)) \in \Big(e^{n^\alpha}(1-2\varepsilon), \ e^{n^\alpha}(1+2\varepsilon)\Big).
\]
Hence,
\[
E_\varphi(\Phi) \subset \bigcup_N B\Big(e^{n^{\alpha/b}}, \ {3\varepsilon \over b} n^{\alpha(1/b-1)}e^{n^{\alpha/b}}, \ N \Big).
\]
On the other hand, for a sequence of positive numbers $\epsilon_n$ converging to $0$, we have
\[
E_\varphi(\Phi) \supset  B\big(e^{n^{\alpha/b}}, \ \varepsilon_n n^{\alpha(1/b-1)}e^{n^{\alpha/b}}, \ 1 \big).
\]
Applying Lemma \ref{lem:A}, we have
\[
\dim_H E_\varphi(\Phi) = \liminf_{\ell\to\infty} \frac {\sum_{j=1}^\ell j^{\alpha/b}}{d\sum_{j=1}^{\ell+1} j^{\alpha/b} - (\ell+1)^{\alpha/b}} = \frac 1d.
\]

\subsection{Theorem \ref{thm-3}, case I-2}
From the formula (\ref{eqn:diff}), we get
\[
\varphi(a_n(x)) \in \Big(e^{n^\alpha}(1-2\varepsilon), \ e^{n^\alpha}(1+2\varepsilon)\Big).
\]
Hence,
\[
E_\varphi(\Phi) \subset \bigcup_N B\Big({n^{\alpha/c}}, \ {3\varepsilon \over c} n^{\alpha(1/c-1)}, \ N \Big).
\]
On the other hand, for a sequence of positive numbers $\epsilon_n$ converging to $0$, we have
\[
E_\varphi(\Phi) \supset  B\big(n^{\alpha/c}, \ \varepsilon_n n^{\alpha(1/c-1)}, \ 1\big).
\]
We then apply Lemma \ref{lem:A} to obtain
\[
\dim_H E_\varphi(\Phi) = \liminf_{\ell\to\infty} \frac {\sum_{j=1}^\ell \alpha (1/c-1) \log j}{d\sum_{j=1}^{\ell+1} \alpha/c \log j - \alpha (1/c-1) \log (\ell+1)} = \frac {1-c} d.
\]

\subsection{Theorem \ref{thm-3}, case II}

From the formula (\ref{eqn:diff}), we get
\[
\varphi(a_n(x)) \in \Big(e^{\beta^n}(1-2\varepsilon), \ e^{\beta^n}(1+2\varepsilon)\Big).
\]
Hence,
\[
E_\varphi(\Phi) \subset \bigcup_N B\Big({\beta^{n/c}}, \  {3\varepsilon \over c} \beta^{n (1/c-1)}, \ N\Big).
\]
On the other hand, for a positive sequence $\epsilon_n$ converging to $0$, we have
\[
E_\varphi(\Phi) \supset  B\big({\beta^{n/c}}, \ \varepsilon_n  \beta^{n (1/c-1)}, \ 1 \big).
\]
Applying Lemma \ref{lem:A}, we obtain
\[
\dim_H E_\varphi(\Phi) = \liminf_{\ell\to\infty} \frac {\sum_{j=1}^\ell j(1/c-1) \log\beta}{d\sum_{j=1}^{\ell+1} j/c \log\beta - (\ell+1)(1/c-1)\log \beta} = \frac {1-c}d.
\]

\subsection{Theorem \ref{thm-3}, case III}

From the formula (\ref{eqn:diff}), we get
\[
\varphi(a_n(x)) \in \Big(e^{e^{\gamma^{n}}}(1-2\varepsilon), \ e^{e^{\gamma^{n}}}(1+2\varepsilon)\Big).
\]
Hence,
\[
E_\varphi(\Phi) \subset \bigcup_N B\Big(e^{{1 \over c}\gamma^{n}}, \ {3\varepsilon \over c} e^{\gamma^{n} (1/c-1)}, \ N\Big).
\]
On the other hand, for a positive sequence $\epsilon_n$ converging to $0$, we have
\[
E_\varphi(\Phi) \supset  B\big(e^{{1 \over c}\gamma^{n}}, \ \varepsilon_n  e^{\gamma^{n} (1/c-1)}, \ 1\big).
\]
Applying Lemma \ref{lem:A}, we get
\[
\dim_H E_\varphi(\Phi) = \liminf_{\ell\to\infty} \frac {\sum_{j=1}^\ell (1/c-1)\gamma^j }{d\sum_{j=1}^{\ell+1} 1/c \gamma^j - (1/c-1)\gamma^{\ell+1}} = \frac {1-c}{d\gamma - (1-c)(\gamma-1)}.
\]

\vskip 5mm

We also apply Lemma \ref{lem:A}  for the lower bounds of Theorem \ref{thm-1}, subcase I-2a and Theorem \ref{thm-2}, subcase I-2a. But for the upper bounds we need Lemma \ref{estimation-upper} and Lemma \ref{estimation-upper-2} respectively.  

\subsection{Proof of Theorem \ref{thm-1}, case I-2a}
We first show the lower bound. Let $x$ be points such that
\[
\varphi(a_n(x)) \in \Big(\alpha n^{\alpha-1}e^{n^\alpha}(1-\varepsilon_n), \ \alpha n^{\alpha-1}e^{n^\alpha}(1+\varepsilon_n)\Big).
\]
where  $(\epsilon_n)_{n\geq 1}$ is a summable positive sequence.
Then 
\[
\sum_{j=1}^n \alpha j^{\alpha-1}e^{j^\alpha}(1-\varepsilon_j) \leq \sum_{j=1}^n\varphi(a_j(x))\leq \sum_{j=1}^n \alpha j^{\alpha-1}e^{j^\alpha}(1+\varepsilon_j),
\]
which implies
\[
e^{n^\alpha}-2 \sum_{j=1}^n \alpha j^{\alpha-1}e^{j^\alpha}\varepsilon_j \leq \sum_{j=1}^n\varphi(a_j(x))\leq e^{n^\alpha}-2 \sum_{j=1}^n \alpha j^{\alpha-1}e^{j^\alpha}\varepsilon_j.
\]
Note that 
\[
\sum_{j=1}^{n/2} \alpha j^{\alpha-1}e^{j^\alpha}\varepsilon_j \leq \sum_{j=1}^{n/2} \alpha j^{\alpha-1}e^{j^\alpha}\leq e^{(n/2)^\alpha},
\]
and by the summability of $(\epsilon_n)_{n\geq 1}$, 
\[
\sum_{j=n/2}^{n} \alpha j^{\alpha-1}e^{j^\alpha}\varepsilon_j \leq \alpha n^{\alpha-1}e^{n^\alpha} \sum_{j=1}^{n/2} \varepsilon_j =o( e^{n^\alpha}).
\]
Hence, these points $x$ are all in $E_\varphi(\Phi)$, that is 
\[
E_\varphi(\Phi) \supset  B\Big((\alpha n^{\alpha-1}e^{n^\alpha})^{1/a}, \ {\varepsilon_n \over a} (\alpha n^{\alpha-1}e^{n^\alpha})^{1/a}, \ 1\Big).
\]
Applying Lemma \ref{lem:A}, we obtain the lower bound.


\bigskip

Now we turn to the upper bound.

Take a subsequence $n_0=1$, and $n_k=\Phi^{-1}(e^k) =k^{1/\alpha}\ (k\geq 1$). If $x\in E_\varphi(\Phi) $
then for any $\epsilon >0$ there exists an integer $N\geq 1$ such that for all $k\geq N$,
\[
(1-\epsilon/5)\Phi(n_k) \leq S_{n_k}\varphi(x) \leq (1+\epsilon/5)\Phi(n_k),
\]
and (as $\Phi(n_k)= e^k$)
\[
(1-\epsilon/5)e^k-(1+\epsilon/5)e^{k-1} \leq S_{n_k}(x) -S_{n_{k-1}}(x) \leq (1+\epsilon/5)e^k-(1-\epsilon/5)e^{k-1}.
\]
Observe that 
\[
(1+\epsilon/5)e^k-(1-\epsilon/5)e^{k-1}< \left((1-\epsilon/5)e^k-(1+\epsilon/5)e^{k-1}\right) \cdot (1+\epsilon).
\]

Fix $\varepsilon=1/3$ and denote by $A_{k}$ the set of points for which the block of symbols $a_{n_{k-1}+1}(x)\cdots a_{n_k}(x)$ in the symbolic expansion of $x$ from the position $n_{k-1}+1$ to $n_{k}$ belongs to the set 
\[
A\left((1-\epsilon/5)e^k-(1+\epsilon/5)e^{k-1}, \ n_k-n_{k-1}, a, \epsilon\right).
\]

Then 
\[
E_\varphi(\Phi) \subset \bigcup_N \bigcap_{k\geq N} A_{ k}.
\]
%
Now, we are going to estimate the upper bound of the Hausdorff dimension of $F=\bigcap_{k\geq 1} A_{ k}$. For $\bigcap_{k\geq N} A_{ k} $ with $N\geq 2$ we have the same bound and the proofs are almost the same.


Let us now define $n(k)=n_k-n_{k-1}$ and $m(k)=(1-\epsilon/5)e^k-(1+\epsilon/5)e^{k-1}$. By the assumption $\alpha>1/2$, we have 
{$
\lim_{k\to \infty} {m(k) \cdot 3^{-n(k)}}=\infty.
$}
Thus we can apply Lemma \ref{estimation-upper} to calculate $G(m(k), n(k),a,1/3, s)$ for all $s>1/d$ and all $k$ large enough. Hence
\begin{align*}
&\sum_{ I_{n_k}(a_1,\ldots,a_{n_k})\cap F\neq\emptyset} |I_{n_k}(a_1,\ldots,a_{n_k})|^{s}\\
 \leq &K_2^{sn_k}\prod_{j=1}^k G(m(j),n(j),a,1/3,s) \\
\leq &  K_2^{sn_k} C_1^k C_2^{n_k-k-1} 3^{-k} \prod_{j=1}^k m(j)^{{1-ds \over a}}.
\end{align*}
As $ds>1$, the right hand side is arbitrarily small for large $k$. This proves
the $s$-dimensional Hausdorff measure
\[
\mathcal{H}^s(F) =0
\]
for all $s>1/d$. We thus obtain the wanted upper bound.

\subsection{Theorem \ref{thm-2}, case I-2a}

For the lower bound, we follow the proof of Theorem \ref{thm-1}, case I-2a, by taking those points $x$ such that
\[
\varphi(a_n(x)) \in \Big(\alpha n^{\alpha-1}e^{n^\alpha}(1-\varepsilon_n), \ \alpha n^{\alpha-1}e^{n^\alpha}(1+\varepsilon_n)\Big).
\]
where  $(\epsilon_n)_{n\geq 1}$ is a summable positive sequence.
Then 
we still have that
these points $x$ are all in $E_\varphi(\Phi)$. By applying the inverse of $\varphi$, we have  
\[
E_\varphi(\Phi) \supset  B\Big(e^{(n^\alpha+\log \alpha +({\alpha-1})\log n)^{1/b}}, \ {2\varepsilon_n \over b} n^{\alpha (1/b-1)}e^{(n^\alpha+\log \alpha +({\alpha-1})\log n)^{1/b}}, \ 1\Big).
\]
Applying Lemma \ref{lem:A}, we obtain the lower bound.

\bigskip

The proof of the upper bound is also similar to that of Theorem \ref{thm-1}, case I-2a. The difference is that we need to apply Lemma \ref{estimation-upper-2} in place of Lemma \ref{estimation-upper}.

As in the proof of Theorem \ref{thm-1}, case I-2a, we take a subsequence $n_0=1$, and $n_k=\Phi^{-1}(e^k) =k^{1/\alpha}\ (k\geq 1)$. Denote by $\widehat{A}_{k}$ the set of points for which the block of symbols $a_{n_{k-1}+1}(x)\cdots a_{n_k}(x)$ in the symbolic expansion of $x$ from the position $n_{k-1}+1$ to $n_{k}$ belongs to the set 
\[
\widehat{A}\left(m(k), n(k), b, 1/3\right),
\]
with $n(k)=n_k-n_{k-1}$ and $m(k)={14\over 15} e^k-{16\over 15}e^{k-1}$.
Then 
\[
E_\varphi(\Phi) \subset \bigcup_N \bigcap_{k\geq N} \widehat{A}_{ k}.
\]

We need only to estimate the upper bound of the Hausdorff dimension of $\widehat F=\bigcap_{k\geq 1} \widehat A_{ k}$.
By the assumption $\alpha>{b \over b+1}>{1 \over 2}$, we still have {$
\lim_{k\to \infty} {m(k) \cdot 3^{-n(k)}}=\infty.
$} 
Thus we can apply Lemma \ref{estimation-upper-2} to calculate $\widehat{G}(m(k), n(k),b,1/3, s)$ for all $s>1/d$ and all $k$ large enough. Hence
\begin{align}
&\sum_{ I_{n_k}(a_1,\ldots,a_{n_k})\cap F\neq\emptyset} |I_{n_k}(a_1,\ldots,a_{n_k})|^{s} \nonumber \\
 \leq &K_2^{sn_k}\prod_{j=1}^k \widehat{G}(m(j),n(j),b,1/3,s) \nonumber  \\
\leq &  K_2^{sn_k} \cdot 6^k \cdot \widehat{C}^{n_k-k-1} \cdot  3^{-k} \prod_{j=1}^k e^{(1-ds)(\log m(j))^{1/b}}. \label{whole-product}
\end{align}
Note that $\log m(j) \approx j$ and $n_k = k^{1/\alpha}$. Thus
{\[
\prod_{j=1}^k e^{(1-ds)(\log m(j))^{1/b}} \stackrel{{e}}{\sim}  e^{(1-ds){b \over b+1}k^{b+1 \over b}} 
\]
 (here we used the symbol $\stackrel{{e}}{\sim}$ appeared in the proof Lemma \ref{lem:A})}. 
As $\frac {b+1} b > \frac 1\alpha$, this is the dominating term of \eqref{whole-product}. As $ds>1$, this term, and the whole product \eqref{whole-product}, converge to $0$ for $k\to\infty$. This proves
the $s$-dimensional Hausdorff measure
\[
\mathcal{H}^s(\widehat F) =0
\]
for all $s>1/d$. We are done.


\subsection*{Acknowledgements} {The authors thank the referee for her/his valuable suggestions which improve the presentation of the paper.} 

The work was supported by the France-Poland bilateral project PHC Polonium (44851YC, PPN/BFR/2019/1/00013). M. R. was also partially supported by National Science Centre grant 2019/33/B/ST1/00275 (Poland).

\bigskip

\bibliographystyle{alpha}

%
%

\end{document}